\documentclass[10pt]{article}
\usepackage{amsmath}
\usepackage{amssymb}
\usepackage{amsthm}
\usepackage{mathrsfs}
\title{Global well-posedness of the 3D Patlak-Keller-Segel system near a straight line}
\author{Bowei Tu \thanks{Email: tbwustc@mail.ustc.edu.cn}}
\date{\today}
\usepackage[left=3cm,right=3cm,top=3cm,bottom=2.5cm]{geometry}
\par 
\begin{document}
	\maketitle 
	\textbf{Abstract.} 
     We consider the Cauchy problem of the three-dimensional parabolic-elliptic Patlak-Keller-Segel chemotactic model. The initial data is almost a Dirac measure supported on a straight line with mass less than $8\pi$. We prove that if the data is sufficiently close to the straight line, then global well-posedness holds.  This result is parallel to the work on vortex filament solutions of the Navier-Stokes equations by Bedrossian, Germain and Harrop-Griffiths \cite{filament}.
	
	\section{Introduction}
	   We consider the following parabolic-elliptic Patlak-Keller-Segel chemotactic system in $\mathbb{R}^d$:
	   \begin{equation}
	   	\begin{cases}
	   		\partial_tu+\nabla\cdot\left(u\nabla c\right)=\Delta u,\\
	   		-\Delta c=u,\\
	   		u|_{t=0}=u_0\ge 0.
	   		
	   	\end{cases}
	   	\label{PKS}
	   \end{equation}
	 where $u$ is cell density and $c$ denotes concentration of chemo-attractant. This system is generally considered the fundamental model for the study of aggregation by chemotaxis of certain microorganisms \cite{Hillen,[45],KS-mold,[49],[59]}. The first equation describes the motion of the microorganism along the gradient of the chemo-attractant. The second equation describes the production and diffusion of the chemo-attractant. From now on, we will refer to (\ref{PKS}) as Patlak-Keller-Segel (PKS).
	 
	 An important property of PKS in two dimensions is that it is $L^1$-critical. That is, if $u(t,x)$ is a solution to (\ref{PKS}), then for all $\lambda\in(0,\infty)$,
	 \begin{equation}
	 	u_{\lambda}(t,x)=\frac{1}{\lambda^2}u\left(\frac{t}{\lambda^2},\frac{x}{\lambda}\right) \label{L^1-scaling}
	 \end{equation}
	 also solves (\ref{PKS}) and preserves the $L^1$ norm. Since the first equation is of divergence form, the mass
	 \begin{equation*}
	 	M:=\int_{\mathbb{R}^2}u_0dx=\int_{\mathbb{R}^2}u(x,t)dx
	 \end{equation*}
	 is a conserved quantity. In fact, the mass serves as a dichotomy between global existence and finite-time blow-up: if $M<8\pi$, then there exists a global in time weak solution \cite{[13],[30]}; if $M>8\pi$, then all solutions blow up in finite time \cite{[13],[46],[54]}; the case $M=8\pi$ is delicate and global well-posedness results for radially symmetric solutions are obtained in \cite{[7]}. We also mention the work by Bedrossian and Masmoudi \cite{2D PKS} that proves local well-posedness in the space of finite Borel measures and establishes significant linear estimates for 2D PKS. There is a remarkable free energy structure for (\ref{PKS}), but we will not discuss it here. Interested readers may refer to \cite{[8],[9],[13]} and references therein for its origin, properties and applications.
	 
	 Another important property of 2D PKS is the existence of self-similar solutions for mass $\alpha\in(0,8\pi)$. These are known to be global attractors for PKS \cite{[13]}. We switch to self-similarity variables:
	 \begin{equation*}
	 	\xi=\frac{x}{\sqrt{t}},\ \ \ \tau=\log t.
	 \end{equation*}
	 and $w(\tau,\xi):=tu(x,t)$. In these variables, (\ref{PKS}) becomes the following
	 \begin{equation}
	 	\begin{cases}
	 		\partial_{\tau}w+\nabla\cdot(w\nabla c)=\Delta w+\frac{1}{2}\nabla\cdot(\xi w),\\
	 		-\Delta c=w.
	 	\end{cases}
	 	\label{self-similar eq}
	 \end{equation}
	 
	  Biler showed in \cite{[5]} that there exists a radially symmetric self-similar solution if and only if the mass $\alpha\in[0,8\pi)$. Naito and Suzuki obtained in \cite{[56]} that the self-similar solution must be radially symmetric by the method of moving planes. Uniqueness is obtained by Biler et al. in \cite{[7]}. We list important properties of self-similar solutions below, which show that in many ways they are qualitatively similar to Oseen vortices of the Navier-Stokes equations. One can see Appendix B in \cite{[13]} for proof of the asymptotic expansion and other interesting properties.
	 
	 \bigskip
	 
	 \textbf{Proposition 1.1.}  Let $\alpha\in(0,8\pi)$.
	 
	 (i) There exists a unique radially symmetric stationary solution to (\ref{self-similar eq}), denoted by $G_{\alpha}$, which is smooth, strictly positive, satisfies $||G_{\alpha}||_1=\alpha$ and denoting $c_{\alpha}=-\Delta^{-1}G_{\alpha}$ we have
	 	\begin{align*}
	 	G_{\alpha}(\xi)&\sim\dfrac{\alpha}{\int e^{c_{\alpha}(\zeta)-|\zeta|^2/4}d\zeta}|\xi|^{-\frac{\alpha}{2\pi}}e^{-|\xi|^2/4}\ \ \ \text{as}\ \xi\rightarrow\infty,\\
	 	\nabla G_{\alpha}(\xi)&\sim-\dfrac{\alpha}{2\int e^{c_{\alpha}(\zeta)-|\zeta|^2/4}d\zeta}\left(\frac{\alpha\xi}{\pi|\xi|^2}+\xi\right)|\xi|^{-\frac{\alpha}{2\pi}}e^{-|\xi|^2/4}\ \ \ \text{as}\ \xi\rightarrow\infty.
	 \end{align*}
	 in the sense of asymptotic expansion. Moreover, for all $p\in[1,\infty]$ and $m\ge 0$ we have
	 	\begin{equation*}
	 	||G_{\alpha}||_{L^p(m)}+||\nabla G_{\alpha}||_{L^p(m)}\lesssim_{m,p,\alpha}1.
	 \end{equation*}
	 
	 (ii) In physical variables, $\frac{1}{t}G_{\alpha}\left(\frac{x}{\sqrt{t}}\right)$ is the unique mild solution with initial data $\alpha\delta$, where $\delta$ denotes the Dirac delta mass. In particular, $\frac{1}{t}G_{\alpha}\left(\frac{x}{\sqrt{t}}\right)\rightharpoonup^*\alpha\delta$ as $t\searrow 0$. 
	 
	 \bigskip
	 
	 When studied in higher dimensions $d\ge 3$, the dynamics for PKS are quite different. The system becomes mass-supercritical and the scaling transformation in (\ref{L^1-scaling}) preserves the $L^{d/2}$ norm. The authors in \cite{3D PKS wellposedness} showed that for initial data $||u_0||_{L^{d/2}}<C(d)$, where $C(d)$ is related to Gagliardo-Nirenberg inequality, there exists a global in time weak solution. See also references in \cite{3D PKS wellposedness} for earlier results and technique concerning the global existence for (\ref{PKS}). For certain size of initial data, there exists finite-time blow-up solutions of (\ref{PKS}), see for example \cite{[26]}. It is worth noting that there still exists a gap between the available criteria ensuring blow-up and global existence for general solutions of PKS. Partial results about radially symmetric case in terms of Morrey norms are recently obtained in \cite{Global radial Morrey,variant diffusion Morrey,3D blow-up,Morrey blow-up}.
	 
	 \bigskip
	 
	 In this paper, we want to understand well-posedness theory for 3D PKS from another viewpoint. Taking $\mathbb{R}^3=\mathbb{R}^2\times\mathbb{R}$ with coordinates $(x,z)\in\mathbb{R}^2\times\mathbb{R}$, an explicit example of solution evolved from the straight line is given by the 2D self-similar solution $\frac{1}{t}G_{\alpha}\left(\frac{x}{\sqrt{t}}\right)$ with mass less than $8\pi$. Our work, to our knowledge, is the first to prove global well-posedness for 3D PKS near large self-similar solutions. To obtain this result, we apply the fluid dynamics regime to the chemotaxis model. This connection may be partially supported by the similarity of works between \cite{2D PKS} and \cite{2D NS uniqueness}. In addition, recently \cite{filament} analyzed the dynamics near the vortex filament in 3D Navier-Stokes regime and we expect similar results can hold for 3D PKS model. However, there is a new difficulty in our work: the perturbation term for core part in NSE case has zero integral. This no longer holds for PKS,  hence the exponential decay estimate for 2D semigroups cannot be directly applied to control the perturbation term. In order to overcome this difficulty, we have to further decompose the perturbation term into two parts: one is of divergence form and the other has fewer nonlocal operators. By controlling the two parts separately, we can obtain the exponential-decay linear estimate for the core part, which is essential for the construction of contraction mapping.
	 
	 \bigskip
	 
	 It is natural to expect that the 2D self-similar solution provides the microscopic structure for the evolution of \textit{any} smooth, non-self-interacting curve. For 3D NSE, this expectation is, in some sense, confirmed by \cite{filament}, where local coordinates are introduced near the filament that straighten out the curve and local well-posedness holds even for large perturbations of filament. For 3D PKS, however, this is still an open problem, and we hope to establish similar results in future papers.
	 
	 \bigskip
	 
	 \textbf{Notations and Conventions}
	 
	 \bigskip
	 
	 We use the following conventions regarding vector calculus:
	 
	 \begin{itemize}
	 	\item  $\left(a_{ij}\right)_{ij}, \left(a_j^i\right)_{ij}, \left(a^{ij}\right)_{ij}$ all denote the matrix with line index $i$, column index $j$.
	 	
	 	\item If $a,b$ are vectors, then $a\otimes b=\left(a^ib^j\right)_{ij}$.
	 	
	 	\item For vectors, we denote $|v|$ to be the usual norm induced by the Euclidean metric.
	 \end{itemize}
	 
	 We denote $g\lesssim f$ if there exists a constant $C>0$ such that $g\le Cf$ and we use $g\lesssim_{\alpha,\beta,...}f$ to emphasize dependence of $C$ on parameters $\alpha,\beta,...$. We similarly write $g\approx f$ if we have both $g\lesssim f$ and $f\lesssim g$. We use $\left\langle\xi\right\rangle=\left(1+|\xi|^2\right)^{1/2}$. Finally, throughout the paper, we only take Fourier transform in $z$.
	 
	\section{Statement of Results and Outline of the Proof}
    \subsection{Function Spaces}
    In order to state our results, we first define several function spaces.
    
    To handle the self-similar part of the solution, for $1\le p<\infty$ and $m\ge 0$ we define the weighted Lebesgue space $L_{\xi}^p(m)$ with norm 
    \begin{equation*}
    	||f||_{L_{\xi}^p(m)}^p=\int_{\mathbb{R}^2}\left\langle\xi\right\rangle^{pm}|f(\xi)|^pd\xi.
    \end{equation*}
    In order to get exponential decay estimates for semigroups, we also introduce the following function spaces:
    \begin{equation*}
    	L_0^2(m)=\left\{f\in L_{\xi}^2(m)\Bigg|\int_{\mathbb{R}^2}fd\xi=0\right\},
    \end{equation*}
    
    The following normalization is used for the Fourier transform in the $z$-direction throughout the paper:
    \begin{equation*}
    	\hat{f}(\zeta)=\frac{1}{\sqrt{2\pi}}\int f(z)e^{-iz\zeta}dz.
    \end{equation*}
    To control the regularity in the translation-invariant $z$-direction, for a Banach space X of functions defined on $\mathbb{R}^2$, we define $B_zX$ as the space of functions defined on $\mathbb{R}^2\times\mathbb{R}$ with norm
    \begin{equation*}
    	||f||_{B_zX}=\int_{\mathbb{R}}||\hat{f}(\cdot,\zeta)||_Xd\zeta,\ \ (x,z)\in\mathbb{R}^2\times\mathbb{R}.
    \end{equation*}
    
    Finally, given a space of functions $X$ continuously embedded in the space of tempered distributions $\mathscr{S}'$, we say that $u\in C_w([0,T],X)$ if $u\in L^{\infty}([0,T],X)$ and for all $t_0\in[0,T]$ and test functions $\phi\in\mathscr{S}$ we have
    \begin{equation*}
    	\lim\limits_{t\rightarrow t_0}\left\langle u(t),\phi\right\rangle=\left\langle u(t_0),\phi\right\rangle.
    \end{equation*}
    
    \subsection{Mild Solution}
    Formally, we write the solution of (\ref{PKS}) with initial data $u(t=0)=u_0$ using the Duhamel formula as
    \begin{equation}
    	u(t)=e^{t\Delta}u_0-\int_{0}^{t}e^{(t-s)\Delta}\text{div}(u(s)\nabla c(s))ds. \label{duhamel formula}
    \end{equation}
    Giga and Miyakawa \cite{mild solution morrey} considered solution of 3D NSE with vorticity in Morrey type space defined as the set of signed measures satisfying $\sup_{r>0,x\in\mathbb{R}^3}|\mu|(B(x,r))/r<+\infty$, where $|\mu|$ is the total variation of $\mu$, and they observed that vortex rings, i,e., measures supported on a closed curve, fall into precisely this class. In fact, this setting can be applied to the chemotaxis model and we give a rigorous definition of what we mean by a mild solution of (\ref{PKS}):
    
    \bigskip
    
    \textbf{Definition 2.1.} Let $M^{\frac{3}{2}}$ be the space of regular Borel measures such that 
    \begin{equation*}
    	||\mu||_{M^{\frac{3}{2}}}:=\sup_{r>0,y\in\mathbb{R}^3}\left\{r^{-1}|\mu(B(y,r))|\right\}<\infty.
    \end{equation*}
    Given a $T>0$, we call a function $u\in C_w([0,T];M^{\frac{3}{2}})$ a mild solution to (\ref{PKS}) with initial data $u_0\in M^{\frac{3}{2}}$ provided
    
    (i) the initial data is attained $u(0)=u_0$ (hence $u(t)\rightharpoonup^*u(0)$ as $t\searrow 0$);
    
    (ii) the equations are satisfied in the sense of Duhamel's formula (\ref{duhamel formula}) (and in particular, the Duhamel integral is well-defined).
    
    \subsection{Main Results}
    Now we state our main results for small perturbations of the straight line:
    
    \bigskip
    
    \textbf{Theorem 2.2.} For any $\alpha\in(0,8\pi)$ and any $m>2$, there exists $\epsilon>0$ such that if $\mu^b:\mathbb{R}^2\times\mathbb{R}\rightarrow\mathbb{R}$ satisfies $||\mu^b||_{B_zL_x^1}\le\epsilon$, then the following holds:
    
    (i) (Existence) There exists a global mild solution of (\ref{PKS}) with initial data
    	\begin{equation}
    	u(t=0)=\alpha\delta_{x=0}+\mu^b \label{main thm initial data}.
    \end{equation}
    which can be decomposed into
    \begin{equation*}
    	u(t,x,z)=\frac{1}{t}G_{\alpha}\left(\frac{x}{\sqrt{t}}\right)+\frac{1}{t}U^c\left(\log t,\frac{x}{\sqrt{t}},z\right)+u^b(t,x,z),
    \end{equation*}
    where the core part $U^c$ and the background part $u^b$ satisfy the estimates
     	\begin{equation}
     	\sup_{-\infty<\tau<\infty}||U^c(\tau)||_{B_zL_\xi^2(m)}+\sup_{0<t<\infty}t^{\frac{1}{4}}||u^b(t)||_{B_zL_x^{4/3}}\lesssim\epsilon. \label{main thm norm}
     \end{equation}
     
     (ii) (Uniqueness) If $u'$ is another mild solution with initial data (\ref{main thm initial data}) admitting the decomposition
     \begin{equation*}
     	u'(t,x,z)=\frac{1}{t}G_{\alpha}\left(\frac{x}{\sqrt{t}}\right)+\frac{1}{t}(U^c)'\left(\log t,\frac{x}{\sqrt{t}},z\right)+(u^b)'(t,x,z),
     \end{equation*}
     where $(U^c)'$ and $(u^b)'$ satisfy the bounds (\ref{main thm norm}), then $u=u'$.
     
     \bigskip
     
     (iii) (Lipschitz dependence) The solution map from the data to solution
     \begin{equation*}
     	\mu^b\mapsto(u^b,U^c)
     \end{equation*}
     is Lipschitz continuous if one endows the data space with norm on $B_zL_x^1$ and the solution space with the norm (\ref{main thm norm}).
     
     \bigskip
     
     The proof of Theorem 2.2 follows from applying the contraction mapping principle to the equations satisfied by the core and background pieces. In order to obtain bounds for these pieces we first introduce the self-similar coordinates
     \begin{equation*}
     	\tau=\log t,\ \ \ \xi=\frac{x}{\sqrt{t}},\ \ \ z=z.
     \end{equation*}
     where we note that as $G_{\alpha}$ is translation-invariant in $z$, we do not rescale the $z$-coordinate. We then define
     \begin{equation*}
     	U(\tau,\xi,z)=e^{\tau}u(e^{\tau},e^{\frac{\tau}{2}}\xi,z),\ \ \ V(\tau,\xi,z)=e^{\frac{\tau}{2}}v(e^{\tau},e^{\frac{\tau}{2}}\xi,z).
     \end{equation*}
     where $v=\nabla(-\Delta)^{-1}u$. We may write the equation (\ref{PKS}) as
     \begin{equation*}
     	\begin{cases}
     		\partial_{\tau}U+\overline{\nabla}\cdot(U V)=(\mathcal{L}+e^{\tau}\partial_{zz})U,\\
     		-\overline{\Delta}c=U.
     	\end{cases}
     \end{equation*}
     where the rescaled gradient $\overline{\nabla}$, rescaled Laplacian $\overline{\Delta}$, and the 2D Fokker-Planck operator $\mathcal{L}$ are defined by
     \begin{equation*}
     	\overline{\nabla}=\begin{bmatrix}
     		\nabla_{\xi}\\
     		e^{\frac{\tau}{2}}\partial_{z}
     	\end{bmatrix},\ \ \ \overline{\Delta}=\Delta_{\xi}+e^{\tau}\partial_{zz},\ \ \ \mathcal{L}=\Delta_{\xi}+\frac{1}{2}\xi\cdot\nabla_{\xi}+1.\ \ \ 
     \end{equation*}
     Finally denote
     \begin{equation*}
     	U^g=G_{\alpha}(\xi),\ \ \ V^g=\begin{bmatrix}
     		V^{G_{\alpha}}(\xi)\\
     		0
     	\end{bmatrix}.
     \end{equation*}
     where $V^{G_{\alpha}}=\nabla_{\xi}(\Delta_{\xi})^{-1}G_{\alpha}$. 
     
     The core piece, $U^c$, is taken to satisfy the equation
     \begin{equation*}
     	\begin{cases}
     		\partial_{\tau}U^c+\overline{\nabla}\cdot((U^g+U^c)\cdot V)=(\mathcal{L}+e^{\tau}\partial_{zz})(U^g+U^c),\\
     		\lim\limits_{\tau\rightarrow-\infty}U^c(\tau)=0.
     	\end{cases}
     \end{equation*}
     In order to solve this equation, we need to prove estimates for the solution operator $U(\tau)=S(\tau,\sigma)U(\sigma)$ of the corresponding linearized equation:
     \begin{equation*}
     	\begin{cases}
     		\partial_{\tau}U+\nabla_{\xi}\cdot(U V^{G_{\alpha}})+\overline{\nabla}\cdot(G_{\alpha}V)=(\mathcal{L}+e^{\tau}\partial_{zz})U,\\
     		V=\overline{\nabla}\left(-\overline{\Delta}\right)^{-1}U.
     	\end{cases}
     \end{equation*}
     
     The background piece, $u^b$, is taken to satisfy the equation
     \begin{equation*}
     	\begin{cases}
     		\partial_{t}u^b+\nabla\cdot(u^bv)=\Delta u^b,\\
     		u^b|_{t=0}=\mu^b.
     	\end{cases}
     \end{equation*}
     Solutions are constructed by establishing estimates for the solution operator $u(t)=\mathbb{S}(t,s)u(s)$ for the corresponding linearized equation
     \begin{equation*}
     	\partial_{t}u+\nabla\cdot(uv^g)=\Delta u
     \end{equation*}
     where $v^g$ is $V^g$ in physical coordinates. 
     
     \bigskip
     
     The bulk of this paper is to obtain estimates for the linear propagators $S(\tau,\sigma),\mathbb{S}(t,s)$. Given these bounds, the proof of Theorem 2.1 follows from an elementary application of the contraction mapping principle that we carry out in Section 5.
     
    \section{Linear Estimates for the Core Part}
    \subsection{Statement of the Estimates}
    Consider the following linearized equation
    \begin{equation}
    	\partial_{\tau}U+\nabla_{\xi}\cdot(U V^{G_{\alpha}})+\overline{\nabla}\cdot(G_{\alpha}V)=(\mathcal{L}+e^{\tau}\partial_{zz})U. \label{core linear}
    \end{equation}
    For $\tau\ge\sigma$ we define the solution operator $S(\tau,\sigma)$ for the equation (\ref{core linear}) by 
    \begin{equation*}
    	U(\tau)=S(\tau,\sigma)U(\sigma).
    \end{equation*}
    
    In this section we prove the following result:
    
    \bigskip
    
    \textbf{Theorem 3.1.} Let $\alpha\in(0,8\pi)$ and $m>2$. Then for all $\sigma\in\mathbb{R}$ the map $\tau\mapsto S(\tau,\sigma)$ is continuous as a map from $[\sigma,\infty)$ to the space of bounded operators on $B_zL_{\xi}^2(m)$. For any scalar $F$ we have the estimate
    \begin{equation}
    	||S(\tau,\sigma)F||_{B_zL_{\xi}^2(m)}\lesssim||F||_{B_zL_{\xi}^2(m)}. \label{core estimate 1}
    \end{equation}
    
    Further, there exists $\nu=\nu(\alpha)\in(0,\frac{1}{2}),$ such that, if $\int_{\mathbb{R}^2} Fd\xi=0$, we have the estimate
    \begin{equation}
    	||S(\tau,\sigma)F||_{B_zL_{\xi}^2(m)}\lesssim e^{-\nu(\tau-\sigma)}||F||_{B_zL_{\xi}^2(m)}. \label{decay core estimate}
    \end{equation}
    
    Finally, for $1<p\le2$ and all 3-vectors $F\in B_zL_{\xi}^p(m)$, we have the estimate
    \begin{equation}
    	||S(\tau,\sigma)\overline{div}F||_{B_zL_{\xi}^2(m)}\lesssim\dfrac{e^{-\nu(\tau-\sigma)}}{a(\tau-\sigma)^{\frac{1}{p}}}||F||_{B_zL_{\xi}^p(m)}. \label{core estimate 2}
    \end{equation}
    where $a(\tau)=1-e^{-\tau}$, $\overline{div}F=\overline{\nabla}\cdot F$.
    
    \bigskip
    
    \subsection{Long Time Estimates}
    In this subsection we prove that the solution operator $S(\tau,\sigma)$ is well-defined, and satisfies the estimates (\ref{core estimate 1}), (\ref{decay core estimate}).
    
    We start by taking the Fourier transform in $z$ of the equation (\ref{core linear}) and setting $W(\tau,\xi,\zeta)=\widehat{U}(\tau,\xi,\zeta)$ to obtain the equation
    \begin{equation}
    	(\partial_{\tau}+e^{\tau}|\zeta|^2-\mathcal{L}+\Lambda_{\alpha})W=\tilde{Z}(W), \label{core Fourier}
    \end{equation}
    where the linear operator is denoted
    \begin{equation*}
    	\Lambda_{\alpha}=\nabla_{\xi}\cdot\left(V^{G_{\alpha}}\cdot\right)+\nabla_{\xi}\cdot V^{G_{\alpha}}+\nabla_{\xi}G_{\alpha}\cdot\nabla_{\xi}\left(\left(-\Delta_{\xi}\right)^{-1}\cdot\right),
    \end{equation*}
    and the perturbation term is given by
    \begin{equation*}
    	Z(W)=-\nabla_{\xi}G_{\alpha}\cdot\nabla_{\xi}\left(\left(-\Delta_{\xi}+e^{\tau}|\zeta|^2\right)^{-1}W\right)+\nabla_{\xi}G_{\alpha}\cdot\nabla_{\xi}\left(\left(-\Delta_{\xi}\right)^{-1}W\right).
    \end{equation*}
    
    The existence of the solution operator $S(\tau,\sigma)$ and the estimates (\ref{core estimate 1}), (\ref{decay core estimate}) are given in the following proposition:
    
    \bigskip
    
    \textbf{Proposition 3.2.} Let $m>2$, $\alpha\in(0,8\pi)$, and $\zeta\in\mathbb{R}$ be fixed. Then, for all $\sigma\in\mathbb{R}$ and all $w_{\sigma}\in L_{\xi}^2(m)$, there exists a unique mild solution $W\in C\left([\sigma,\infty);L_{\xi}^2(m)\right)$ satisfying $W(\sigma)=w_{\sigma}$ and the estimate
    \begin{equation}
    	||W(\tau)||_{L_{\xi}^2(m)}\lesssim||w_{\sigma}||_{L_{\xi}^2(m)} \label{prop 3.2 i}.
    \end{equation}
    
    Furthermore, there exists some $\nu=\nu(\alpha)\in\left(0,\frac{1}{2}\right)$ so that, for all $w_{\sigma}\in L_0^2(m)$, we have the improved estimate
    \begin{equation}
    	||W(\tau)||_{L_{\xi}^2(m)}\lesssim e^{-\nu(\tau-\sigma)}||w_{\sigma}||_{L_{\xi}^2(m)} \label{prop 3.2 ii}.
    \end{equation}
    In both estimates, the implicit constant is independent of $\zeta$.
    
    \bigskip
    
    We first recall some known properties of the semigroup $e^{\tau\mathcal{L}}$:
    
    \bigskip
    
    \textbf{Lemma 3.3} \cite{Fokker-Planck}. Let $m\ge 0$. 
    
    (i) $e^{\tau\mathcal{L}}$ is a strongly continuous semigroup on $L^2(m)$ and satisfies
    \begin{equation}
    	||e^{\tau\mathcal{L}}f||_{L^2(m)}\lesssim||f||_{L^2(m)},\ \ \ 
    	||\nabla e^{\tau\mathcal{L}}f||_{L^2(m)}\lesssim\dfrac{1}{a(\tau)^{\frac{1}{2}}}||f||_{L^2(m)}.
    \end{equation}
    
    (ii) If $m>2$, then for all $f\in L_0^2(m)$,
    \begin{equation}
    	||e^{\tau\mathcal{L}}f||_{L^2(m)}\lesssim e^{-\frac{1}{2}\tau}||f||_{L^2(m)}.
    \end{equation}
    
    (iii) Furthermore, for all $1\le q\le p \le \infty$ and $\alpha\in\mathbb{N}^2$,
    \begin{equation}
    	||\partial^{\alpha}e^{\tau\mathcal{L}}f||_{L^p(m)}\lesssim\dfrac{1}{a(\tau)^{\frac{1}{q}-\frac{1}{p}+\frac{|\alpha|}{2}}}||f||_{L^q(m)}.
    \end{equation}
    
    As a compact perturbation of $e^{\tau\mathcal{L}}$, $e^{\tau(\mathcal{L}-\Lambda_{\alpha})}$ has similar semigroup estimates:
    
    \bigskip
    
    \textbf{Lemma 3.4} \cite{2D PKS}. Let $\alpha\in(0,8\pi)$ and $m>2$.
    
    (i) $e^{\tau(\mathcal{L}-\Lambda_{\alpha})}$ is a strongly continuous semigroup on $L_{\xi}^2(m)$ and satisfies
   \begin{equation}
   	||e^{\tau(\mathcal{L}-\Lambda_{\alpha})}f||_{L^2(m)}\lesssim_{\alpha}||f||_{L^2(m)},\ \ \ ||\nabla e^{\tau(\mathcal{L}-\Lambda_{\alpha})}f||_{L^2(m)}\lesssim_{\alpha}\dfrac{1}{a(\tau)^{\frac{1}{2}}}||f||_{L^2(m)}  \label{2D core 1}
   \end{equation}

    (ii) There exists $\nu=\nu(\alpha)\in(0,1/2)$ such that for all $f\in L_0^2(m)$,
    \begin{equation}
    	||e^{\tau(\mathcal{L}-\Lambda_{\alpha})}f||_{L^2(m)}\lesssim_{\alpha} e^{-\nu\tau}||f||_{L^2(m)}. \label{2D core 2}
    \end{equation}
    
    (iii) If $q\in(1,2]$ then $e^{\tau(\mathcal{L}-\Lambda_{\alpha})}\nabla$ is a bounded operator from $L^q(m)$ to $L_0^2(m)$ and there exists a $\nu\in(0,1/2)$ (the same $\nu$ as in (ii)) such that 
    \begin{equation}
    	||e^{\tau(\mathcal{L}-\Lambda_{\alpha})}\nabla f||_{L^2(m)}\lesssim_{\alpha}\dfrac{e^{-\nu\tau}}{a(\tau)^{1/q}}||f||_{L^q(m)}.
    \end{equation}
    
    \bigskip
    
    In addition, we prove several Biot-Savart type estimates that are crucial for the control of perturbation terms:
    
    \bigskip
    
    \textbf{Lemma 3.5.} Let $m>1$ and $\lambda>0$.
    
    (i) If $1<r\le\infty$ and $0<\delta\le\min\{\frac{1}{2},1-\frac{1}{r}\}$, then 
    \begin{equation}
    	||\left(\lambda^2-\Delta_{\xi}\right)^{-1}f||_{L_{\xi}^r}\lesssim\lambda^{-\left(\frac{2}{r}+2\delta\right)}||f||_{L_{\xi}^2(m)}. 
    \end{equation}
    
    (ii) If $2\le r<\infty$, then
    \begin{equation}
    	||\nabla_{\xi}\left(-\Delta_{\xi}\right)^{-1}f||_{L_{\xi}^r}\lesssim||f||_{L_{\xi}^2(m)}.
    \end{equation}
    
    (iii)      	\begin{equation}
    	||\nabla\left(-\Delta\right)^{-1}f||_{B_zL_x^4}\lesssim||f||_{B_zL_x^{4/3}}  \label{Biot-Savart law}
    \end{equation}
    
    \begin{proof}
    	(i) It is known that
    	\begin{equation*}
    		\left(\lambda^2-\Delta_{\xi}\right)^{-1}f(\xi)=\int_{\mathbb{R}^2}K(\lambda(\xi-\eta))f(\eta)d\eta.
    	\end{equation*}
    	where $K$ is the 2D Bessel potential and for all $1\le p<\infty$, $K\in L_{\xi}^p$. If $\frac{1}{p}+\frac{1}{q}=1+\frac{1}{r}$, then applying Young's inequality, we have
    	\begin{equation*}
    		||\left(\lambda^2-\Delta_{\xi}\right)^{-1}f||_{L_{\xi}^r}\lesssim\lambda^{-\frac{2}{p}}||K||_{L_{\xi}^p}||f||_{L_{\xi}^q}.
    	\end{equation*}
    	
    	For $1\le q\le 2$ and $\frac{1}{q}<\frac{m+1}{2}$ we have the embedding
    		\begin{equation*}
    		||f||_{L_{\xi}^q}\lesssim||f||_{L_{\xi}^2(m)}.
    	\end{equation*}
    	
    	Taking $\frac{1}{p}=\frac{1}{r}+\delta$ and $\frac{1}{q}=1-\delta$, we obtain the estimate in (i).
    	
    	\bigskip
    	
    	(ii) Note that $\nabla_{\xi}\left(-\Delta_{\xi}\right)^{-1}f=k\ast f$, where
    		\begin{equation*}
    		k(\xi)=-\frac{1}{2\pi}\frac{\xi}{|\xi|^2}.
    	\end{equation*}
    	
    	Hardy-Littlewood-Sobolev's inequality tells us that
    	\begin{equation*}
    		||\nabla_{\xi}\left(-\Delta_{\xi}\right)^{-1}f||_{L_{\xi}^r}\lesssim||f||_{L^{\frac{2r}{2+r}}}.
    	\end{equation*}
    	
    	Using the embedding $L_{\xi}^2(m)\subset L^{\frac{2r}{2+r}}$, we obtain the estimate in (ii).
    	
    	\bigskip
    	
    	(iii) By scaling it suffices to prove that
    	\begin{equation*}
    		||\nabla_x(1-\Delta_x)^{-1}f||_{L_x^4}\lesssim||f||_{L_x^{4/3}},\ \ \ \ \ \ ||(1-\Delta_x)^{-1}f||_{L_x^4}\lesssim||f||_{L_x^{4/3}}.
    	\end{equation*}
    	This follows from Young's inequality and fact that the kernels of $\nabla_x(1-\Delta_x)^{-1}$ and $(1-\Delta_x)^{-1}$ are both in $L^2$.
    \end{proof}
    \bigskip
    
    Now we can control the perturbation term $Z(W)$:
    
    \bigskip
    
    \textbf{Proposition 3.6.} If $m>1$, $0<\delta<\frac{1}{2}$, then we have the estimate
     \begin{equation}
     	||Z(W)||_{L_{\xi}^2(m)}\lesssim|\zeta|^{(1-2\delta)}e^{(\frac{1}{2}-\delta)\tau}||W||_{L_{\xi}^2(m)}.
     \end{equation}

     \begin{proof}
     	 Note that
     	 \begin{align*}
     	 	Z(W)&=\nabla_{\xi}G_{\alpha}\cdot\nabla_{\xi}\left(\left(-\Delta_{\xi}\right)^{-1}W-\left(-\Delta_{\xi}+e^{\tau}|\zeta|^2\right)^{-1}W\right)\\
     	 	&=\nabla_{\xi}G_{\alpha}\cdot\nabla_{\xi}\left(e^{\tau}|\zeta|^2\left(-\Delta_{\xi}\right)^{-1}\left(-\Delta_{\xi}+e^{\tau}|\zeta|^2\right)^{-1}W\right)
     	 \end{align*}
      Applying Young's inequality and Lemma 3.5 (ii) successively, we obtain
     	\begin{align*}
     		&\left|\left|\nabla_{\xi}\left(e^{\tau}|\zeta|^2\left(-\Delta_{\xi}\right)^{-1}\left(-\Delta_{\xi}+e^{\tau}|\zeta|^2\right)^{-1}W\right)\right|\right|_{L_{\xi}^{\infty}}\\
     		\le &\left(e^{\tau}|\zeta|^2\right)^{\frac{1}{2}-\delta}\left|\left|\nabla_{\xi}\left(-\Delta_{\xi}\right)^{-1}W\right|\right|_{L_{\xi}^{\frac{2}{1-2\delta}}}\\
     		\lesssim&\left(e^{\tau}|\zeta|^2\right)^{\frac{1}{2}-\delta}||W||_{L_{\xi}^2(m)}.
     	\end{align*}
     	
     		Combining with the fast decay of $\nabla_{\xi}G_{\alpha}$, we finally get the desired estimate.
     \end{proof}
     We proceed to the proof of Proposition 3.2.
     \begin{proof}
     	We first write the integral equation for (\ref{core Fourier}):
     	\begin{equation}
     		W(\tau)=e^{-(e^{\tau}-e^{\sigma})|\zeta|^2}e^{(\tau-\sigma)(\mathcal{L}-\Lambda_{\alpha})}W(\sigma)+\int_{\sigma}^{\tau}e^{-|\zeta|^2(e^{\tau}-e^s)}e^{(\tau-s)(\mathcal{L}-\Lambda_{\alpha})}\tilde{Z}(W(s))ds. \label{core Fourier integral eq}
     	\end{equation}
     	Take $\epsilon>0$ and define the map
     	\begin{equation*}
     		T:C\left([\sigma,\sigma+\epsilon];L_{\xi}^2(m)\right)\rightarrow C\left([\sigma,\sigma+\epsilon];L_{\xi}^2(m)\right),
     	\end{equation*}
     	by
     	\begin{equation*}
     		T(W)=\int_{\sigma}^{\tau}e^{-|\zeta|^2(e^{\tau}-e^s)}e^{(\tau-s)(\mathcal{L}-\Lambda_{\alpha})}\tilde{Z}(W(s))ds.
     	\end{equation*}
     	
     	Taking $\delta=\frac{1}{4}$ and applying the above lemma for the perturbative term $Z$, we may then bound
     		\begin{align*}
     		||T(W)||_{L_{\xi}^2(m)}&\lesssim\int_{\sigma}^{\tau}e^{-|\zeta|^2(e^{\tau}-e^s)}|\zeta|^{\frac{1}{2}}e^{\frac{1}{4}s}||W(s)||_{L_{\xi}^2(m)}ds\\
     		&=\int_{\sigma}^{\tau}e^{-|\zeta|^2(e^{\tau}-e^s)}|\zeta|^{\frac{1}{2}}\left(e^{\tau}-e^s\right)^{\frac{1}{4}}\left(\dfrac{e^s}{e^{\tau}-e^s}\right)^{\frac{1}{4}}||W(s)||_{L_{\xi}^2(m)}ds\\
     		&\lesssim\int_{\sigma}^{\tau}\left(\dfrac{1}{e^{\tau-s}-1}\right)^{\frac{1}{4}}ds\cdot||W||_{C\left([\sigma,\sigma+\epsilon];L_{\xi}^2(m)\right)}\\
     		&\lesssim\left(\tau-\sigma\right)^{\frac{3}{4}}||W||_{C\left([\sigma,\sigma+\epsilon];L_{\xi}^2(m)\right)}.
     	\end{align*}
     	
     	By choosing $0<\epsilon=\epsilon(\alpha)<<1$ sufficiently small (independent of $\sigma$), we may use the contraction principle to find a unique mild solution of (\ref{core Fourier}) on the time interval $[\sigma,\sigma+\epsilon]$. Further, as $\epsilon$ is independent of $\sigma$, we may iterate this argument to obtain a global solution.
     	
     	To obtain the first estimate, we use an identical argument to obtain the following control:
     	\begin{equation*}
     		e^{e^{\tau}|\zeta|^2}||W(\tau)||_{L_{\xi}^2(m)}\lesssim e^{e^{\sigma}|\zeta|^2}||W(\sigma)||_{L_{\xi}^2(m)}+\int_{\sigma}^{\tau}e^{e^s|\zeta|^2}|\zeta|^{\frac{1}{2}}e^{\frac{1}{4}s}||W(s)||_{L_{\xi}^2(m)}ds.
     	\end{equation*}
     	
     	Applying the integrated form of Gronwall's inequality to the continuous non-negative function $\tau\mapsto e^{e^{\tau}|\zeta|^2}||W(\tau)||_{L_{\xi}^2(m)}$ we obtain
     		\begin{equation*}
     		e^{e^{\tau}|\zeta|^2}||W(\tau)||_{L_{\xi}^2(m)}\lesssim e^{e^{\sigma}|\zeta|^2}||W(\sigma)||_{L_{\xi}^2(m)}e^{4C|\zeta|^{\frac{1}{2}} (e^{\frac{1}{4}\tau}-e^{\frac{1}{4}\sigma})}.
     	\end{equation*}
     	from which the first estimate follows.
     	
     	To prove the exponential decay estimate (\ref{decay core estimate}), we need to further decompose the perturbation term
     	\begin{equation*}
     		Z(W)=Z_1+Z_2
     	\end{equation*}
     	where
     	\begin{align*}
     		Z_1&:=\nabla_{\xi}\cdot \left(G_{\alpha}\nabla_{\xi}\left(e^{\tau}|\zeta|^2\left(-\Delta_{\xi}\right)^{-1}\left(-\Delta_{\xi}+e^{\tau}|\zeta|^2\right)^{-1}W\right)\right)\\
     		Z_2&:=G_{\alpha}e^{\tau}|\zeta|^2\left(-\Delta_{\xi}+e^{\tau}|\zeta|^2\right)^{-1}W
     	\end{align*}
     	These terms can be controlled as $Z(W)$ itself:
     	\begin{align*}
     		||Z_2||_{L_{\xi}^2(m)}&\lesssim e^{\tau}|\zeta|^2||(-\Delta_{\xi}+e^{\tau}|\zeta|^2)^{-1}W||_{L_{\xi}^2}\\
     		&\lesssim|\zeta|^{(1-2\delta)}e^{(\frac{1}{2}-\delta)\tau}||W||_{L_{\xi}^2(m)}. 
     	\end{align*}
     	
     	\begin{equation*}
     		||Z_1||_{L_{\xi}^2(m)}\lesssim|\zeta|^{(1-2\delta)}e^{(\frac{1}{2}-\delta)\tau}||W||_{L_{\xi}^2(m)}.
     	\end{equation*}
     	
     	Now suppose $w_{\sigma}\in L_0^2(m)$. Notice that for every $\beta\in\mathbb{N}$, we have
     	\begin{equation*}
     		(e^{\frac{\tau}{2}}|\zeta|)^{\beta}e^{-\frac{1}{2}(e^{\tau}-e^{\sigma})|\zeta|^2}\lesssim a(\tau-\sigma)^{-\frac{\beta}{2}}
     	\end{equation*}
     	Using the control for perturbation terms and the above estimate, we obtain
     	 \begin{equation*}
     		||e^{-(e^{\tau}-e^{\sigma})|\zeta|^2}e^{(\tau-\sigma)(\mathcal{L}-\Lambda_{\alpha})}W(\sigma)||_{L_{\xi}^2(m)}\lesssim e^{-\frac{1}{2}(e^{\tau}-e^{\sigma})|\zeta|^2}e^{-\nu(\tau-\sigma)}||W(\sigma)||_{L_{\xi}^2(m)},
     	\end{equation*}
     	
     	 \begin{equation*}
     		||e^{-|\zeta|^2(e^{\tau}-e^s)}e^{(\tau-s)(\mathcal{L}-\Lambda_{\alpha})}Z_1||_{L_{\xi}^2(m)}\lesssim e^{-\frac{1}{2}(e^{\tau}-e^{s})|\zeta|^2}e^{-\nu(\tau-s)}|\zeta|^{\frac{1}{2}}e^{\frac{1}{4}s}||W(s)||_{L_{\xi}^2(m)},
     	\end{equation*}
     	
     	As for the $Z_2$ part, we deal with different time interval separately:
     	
     	\underline{\textit{Case 1:} $s\le\tau\le s+1$}. Then we directly use the estimate (\ref{2D core 1}):
     	\begin{equation*}
     		||e^{-\frac{1}{2}|\zeta|^2(e^{\tau}-e^s)}e^{(\tau-s)(\mathcal{L}-\Lambda_{\alpha})}
     		Z_2||_{L_{\xi}^2(m)}\lesssim |\zeta|^{\frac{1}{2}}e^{\frac{1}{4}s}||W(s)||_{L_{\xi}^2(m)}
     		\lesssim e^{-\nu(\tau-s)}e^{\frac{1}{4}s}|\zeta|^{\frac{1}{2}}||W(s)||_{L_{\xi}^2(m)}.
     	\end{equation*}
     	
     	\underline{\textit{Case 2:} $s+1<\tau\le\sigma$}. In this case $a(\tau-s)^{-\frac{1}{2}}\lesssim 1$, hence
     	
     	\begin{align*}
     		&||e^{-\frac{1}{2}|\zeta|^2(e^{\tau}-e^s)}e^{(\tau-s)(\mathcal{L}-\Lambda_{\alpha})}
     		Z_2||_{L_{\xi}^2(m)}\\
     		=& ||e^{-\frac{1}{2}(\tau-s)}e^{\frac{1}{2}\tau}|\zeta|e^{-\frac{1}{2}(e^{\tau}-e^s)|\zeta|^2+(\tau-\sigma)(\mathcal{L}-\Lambda_{\alpha})}(G_{\alpha}e^{\frac{1}{2}s}|\zeta|(-\Delta_{\xi}+e^s|\zeta|^2)^{-1}W(s))||_{L_{\xi}^2(m)}\\
     		\lesssim& e^{-\frac{1}{2}(\tau-s)}||G_{\alpha}e^{\frac{1}{2}s}|\zeta|(-\Delta_{\xi}+e^s|\zeta|^2)^{-1}W(s)||_{L_{\xi}^2(m)}\\
     		\lesssim& e^{-\nu(\tau-s)}e^{\frac{1}{2}s}|\zeta|\left|\left|(-\Delta_{\xi}+e^s|\zeta|^2)^{-1}W(s)\right|\right|_{L_{\xi}^8}\\
     		\lesssim& e^{-\nu(\tau-s)}e^{\frac{1}{4}s}|\zeta|^{\frac{1}{2}}||W(s)||_{L_{\xi}^2(m)}.
     	\end{align*}
     	Taking the $L_{\xi}^2(m)$ norm on both sides of the integral equation (\ref{core Fourier integral eq}) and using the above bounds, we have
     	\begin{equation*}
     		e^{\frac{1}{2}e^{\tau}|\zeta|^2}e^{v\tau}||W(\tau)||_{L_{\xi}^2(m)}\lesssim e^{\frac{1}{2}e^{\sigma}|\zeta|^2}e^{\nu\sigma}||W(\sigma)||_{L_{\xi}^2(m)}+2\int_{\sigma}^{\tau}e^{\frac{1}{2}e^s|\zeta|^2}e^{\nu s}|\zeta|^{\frac{1}{2}}e^{\frac{1}{4}s}||W(s)||_{L_{\xi}^2(m)}ds.
     	\end{equation*}
     	Applying the integrated form of Gronwall's inequality again to the function $\tau\mapsto e^{\frac{1}{2}e^{\tau}|\zeta|^2+\nu\tau}||W(\tau)||_{L_{\xi}^2(m)}$ as in the proof for (\ref{core estimate 1}), we obtain exponential decay estimate (\ref{decay core estimate}).
     \end{proof}
     
     \subsection{Short Time Estimates}
     In order to prove the estimate (\ref{core estimate 2}), we will combine the long time estimates (\ref{core estimate 1}), (\ref{decay core estimate}) with several short time smoothing estimates. We start with the following estimate motivated by \cite{2D NS uniqueness}:
     
     \bigskip
     
     \textbf{Proposition 3.7.} Let $1<p\le 2$, $m>\frac{1}{2}$. Then there exists some $0<\delta=\delta(\alpha)<<1$ so that for all $\sigma\le\tau\le\sigma+\delta$ and 3-vector $F$ we have the estimate
     \begin{equation}
     	||S(\tau,\sigma)\overline{div}F||_{B_zL_{\xi}^2(m)}\lesssim\dfrac{1}{\left(\tau-\sigma\right)^{\frac{1}{p}}}||F||_{B_zL_{\xi}^p(m)}. \label{core short time estimate}
     \end{equation}
     
     Further, for $\tau>\sigma$, there exists a bounded operator $R(\tau,\sigma)$ on $B_zL_{\xi}^p(m)$ so that $S(\tau,\sigma)\overline{div}F=\overline{div}R(\tau,\sigma)F$ and we have the estimates
     \begin{align}
     	||R(\tau,\sigma)F||_{B_zL_{\xi}^2(m)}&\lesssim\dfrac{1}{\left(\tau-\sigma\right)^{\frac{1}{p}-\frac{1}{2}}}||F||_{B_zL_{\xi}^p(m)},\\
     	||\overline{\nabla}R(\tau,\sigma)F||_{B_zL_{\xi}^2(m)}&\lesssim\dfrac{1}{\left(\tau-\sigma\right)^{\frac{1}{p}}}||F||_{B_zL_{\xi}^p(m)}.
     \end{align}
     
     \begin{proof}
     	Start with the following equation:
     	\begin{equation}
     		\partial_{\tau}F-\left(\mathcal{L}+e^{\tau}\partial_{zz}-\frac{1}{2}\right)F=RHS(F) \label{core R}
     	\end{equation}
     	where $RHS(F)=-\overline{div}F\cdot\begin{bmatrix}
     		V^{G_{\alpha}}\\
     		0
     	\end{bmatrix}-G_{\alpha}\left(\left(-\overline{\Delta}\right)^{-1}\overline{\nabla}\left(\overline{div}F\right)\right)$.
     	
     	We then take $R(\tau,\sigma)$ to be the solution operator for the equation (\ref{core R}), which satisfies $S(\tau,\sigma)\overline{div}=\overline{div}R(\tau,\sigma)$.
     	
     	The solution of (\ref{core R}) may be written using the Duhamel formula as 
     	\begin{equation*}
     		F(\tau)=e^{\left(\mathcal{L}-\frac{1}{2}\right)(\tau-\sigma)+(e^{\tau}-e^{\sigma})\partial_{zz}}F(\sigma)+\int_{\sigma}^{\tau}e^{\left(\mathcal{L}-\frac{1}{2}\right)(\tau-s)+(e^{\tau}-e^s)\partial_{zz}}RHS(F(s))ds,
     	\end{equation*}
     	Our strategy will now be to apply contraction principle to the mapping
     	\begin{equation*}
     		F\mapsto e^{\left(\mathcal{L}-\frac{1}{2}\right)\left(\tau-\sigma\right)+\left(e^{\tau}-e^{\sigma}\right)\partial_{zz}}F(\sigma)+\int_{\sigma}^{\tau}e^{\left(\mathcal{L}-\frac{1}{2}\right)(\tau-s)+\left(e^{\tau}-e^s\right)\partial_{zz}}RHS(F(s))ds
     	\end{equation*}
     	in the closed subspace $X\subset C\left([\sigma,\sigma+\delta];B_zL_{\xi}^p(m)\right)$ with finite norm
     	\begin{equation*}
     		||F||_{X}=\sup_{\tau\in[\sigma,\sigma+\delta]}\left(||F(\tau)||_{B_zL_{\xi}^p(m)}+\left(\tau-\sigma\right)^{\frac{1}{p}-\frac{1}{2}}||F(\tau)||_{B_zL_{\xi}^2(m)}+\left(\tau-\sigma\right)^{\frac{1}{p}}||\overline{\nabla}F(\tau)||_{B_zL_{\xi}^2(m)}\right).
     	\end{equation*}
     	where $\delta>0$ will be chosen sufficiently small (independent of $\sigma$).
     	
     	Noting, on the one hand, that the operator norm on any $L^p$ space of $(e^{\frac{\tau}{2}}\partial_{z})^{\beta}e^{e^{\tau}a(\tau-\sigma)\partial_{z}^2}$ is $\lesssim a(\tau-\sigma)^{-\frac{\beta}{2}}$ and, on the other hand, using semigroup estimates of $e^{\tau\mathcal{L}}$, for $\beta\in\mathbb{N}^3$ and $1\le p\le q\le\infty$ we obtain
     	\begin{equation}
     		\left|\left|\overline{\nabla}^{\beta}e^{(\tau-\sigma)\mathcal{L}+\left(e^{\tau}-e^{\sigma}\right)\partial_z^2}f\right|\right|_{B_zL_{\xi}^q(m)}\lesssim\dfrac{1}{\left(\tau-\sigma\right)^{\frac{1}{p}-\frac{1}{q}+\frac{|\beta|}{2}}}||f||_{B_zL_{\xi}^p(m)}, \label{core duhamel estimate}
     	\end{equation}
     	(Recall that $|\tau-\sigma|\le\delta<<1$ hence $a(\tau-\sigma)\approx\tau-\sigma$.) In particular,
     		\begin{equation*}
     		\left|\left|e^{(\tau-\sigma)\left(\mathcal{L}-\frac{1}{2}\right)+\left(e^{\tau}-e^{\sigma}\right)\partial_{z}^2}F(\sigma)\right|\right|_{X}\lesssim\left|\left|F(\sigma)\right|\right|_{B_zL_{\xi}^p(m)}.
     	\end{equation*}
     	
     	As a consequence, it remains to show that the map
     	\begin{equation*}
     		\mathcal{T}:F\mapsto\int_{\sigma}^{\tau}e^{\left(\mathcal{L}-\frac{1}{2}\right)(\tau-s)+\left(e^{\tau}-e^s\right)\partial_{zz}}RHS(F(s))ds
     	\end{equation*}
     	is a contraction on X. By Lemma 3.5 (iii),
     		\begin{equation*}
     		\left|\left|\left(-\overline{\Delta}\right)^{-1}\overline{\nabla}\left(\overline{div}F\right)\right|\right|_{B_zL_{\xi}^4}\lesssim\left|\left|\overline{div}F\right|\right|_{B_zL_{\xi}^{4/3}}\lesssim\left|\left|\overline{div}F\right|\right|_{B_zL_{\xi}^2(m)}.
     	\end{equation*}
     	Therefore if $1\le p\le 2$, 
     	\begin{align*}
     		||RHS(F)||_{B_zL_{\xi}^p(m)}&\lesssim\left|\left|\overline{div}F\right|\right|_{B_zL_{\xi}^2(m)}\left|\left|V^{G_{\alpha}}\right|\right|_{L_{\xi}^{\frac{2p}{2-p}}}+\left|\left|\overline{div}F\right|\right|_{B_zL_{\xi}^2(m)}\left|\left|G_{\alpha}\right|\right|_{L_{\xi}^{\frac{4p}{4-p}}}\\
     		&\lesssim\left|\left|\overline{div}F\right|\right|_{B_zL_{\xi}^2(m)}.
     	\end{align*}
     	Using the above control over $RHS(F)$, we obtain
     	\begin{align*}
     		||\mathcal{T}(F)||_{B_zL_{\xi}^p(m)}&\lesssim\int_{\sigma}^{\tau}||RHS(F(s))||_{B_zL_{\xi}^p(m)}ds\\
     		&\lesssim\int_{\sigma}^{\tau}\left|\left|\overline{div}F(s)\right|\right|_{B_zL_{\xi}^2(m)}ds\\
     		&\lesssim||F||_{X}\int_{\sigma}^{\tau}\dfrac{1}{\left(s-\sigma\right)^{\frac{1}{p}}}ds\\
     		&\lesssim\left(\tau-\sigma\right)^{1-\frac{1}{p}}||F||_{X},
     	\end{align*}
     	\begin{align*}
     		\left(\tau-\sigma\right)^{\frac{1}{p}}||\overline{\nabla}\mathcal{T}(F)||_{B_zL_{\xi}^2(m)}&\lesssim\left(\tau-\sigma\right)^{\frac{1}{p}}\int_{\sigma}^{\tau}\dfrac{1}{\left(\tau-s\right)^{\frac{1}{p}}}||RHS(F(s))||_{B_zL_{\xi}^p(m)}ds\\
     		&\lesssim\left(\tau-\sigma\right)^{\frac{1}{p}}||F||_{X}\int_{\sigma}^{\tau}\dfrac{1}{\left(\tau-s\right)^{\frac{1}{p}}\left(s-\sigma\right)^{\frac{1}{p}}}ds\\
     		&\lesssim\left(\tau-\sigma\right)^{1-\frac{1}{p}}||F||_{X},
     	\end{align*}
     	\begin{align*}
     		\left(\tau-\sigma\right)^{\frac{1}{p}-\frac{1}{2}}||\mathcal{T}(F)||_{B_zL_{\xi}^2(m)}&\lesssim\left(\tau-\sigma\right)^{\frac{1}{p}-\frac{1}{2}}\int_{\sigma}^{\tau}\dfrac{1}{\left(\tau-s\right)^{\frac{1}{p}-\frac{1}{2}}}||RHS(F(s))||_{B_zL_{\xi}^p(m)}ds\\
     		&\lesssim\left(\tau-\sigma\right)^{\frac{1}{p}-\frac{1}{2}}||F||_{X}\int_{\sigma}^{\tau}\dfrac{1}{\left(\tau-s\right)^{\frac{1}{p}-\frac{1}{2}}\left(s-\sigma\right)^{\frac{1}{p}}}ds\\
     		&\lesssim\left(\tau-\sigma\right)^{1-\frac{1}{p}}||F||_{X}.
     	\end{align*}
     	Combining these estimates we obtain
     	\begin{equation*}
     		||\mathcal{T}(F)||_{X}\lesssim\delta^{1-\frac{1}{p}}||F||_{X}.
     	\end{equation*}
     	so we may choose $0<\delta=\delta(\alpha)<<1$ sufficiently small (independently of $\sigma$) to ensure that $\mathcal{T}$ is a contraction on X. The estimates are then a consequence of the bounds for the solution $F$. 
     \end{proof}
     An essentially identical argument applied directly to the equation (\ref{core linear}) then yields our second short time smoothing estimate:
     
     \bigskip
     
     \textbf{Proposition 3.8.} Let $\alpha\in(0,8\pi)$, $m>\frac{1}{2}$. Then there exists $0<\delta=\delta(\alpha)<<1$ such that for all $\sigma\le\tau\le\sigma+\delta$ we have the estimate
     	\begin{equation}
     	||\overline{\nabla}S(\tau,\sigma)F||_{B_zL_{\xi}^2(m)}\lesssim\dfrac{1}{\left(\tau-\sigma\right)^{\frac{1}{2}}}||F||_{B_zL_{\xi}^2(m)}.
     \end{equation}
     \begin{proof}
     	We first write the equation (\ref{core linear}) in the form
     	\begin{equation*}
     		\partial_{\tau}\Omega-(\mathcal{L}+e^{\tau}\partial_{zz})\Omega=RHS(\Omega),
     	\end{equation*}
     	where
     	\begin{equation*}
     		RHS(\Omega)=-\nabla_{\xi}\cdot\left(\Omega V^{G_{\alpha}}\right)-\overline{\nabla}\cdot\left(G_{\alpha}V\right),
     	\end{equation*}
     	Following a similar argument to Proposition 3.7 we will solve this by applying contraction principle to the mapping
     	\begin{equation*}
     		\Omega\mapsto e^{(\tau-\sigma)\mathcal{L}+\left(e^{\tau}-e^{\sigma}\right)\partial_{zz}}\Omega(\sigma)+\int_{\sigma}^{\tau}e^{(\tau-s)\mathcal{L}+\left(e^{\tau}-e^s\right)\partial_{zz}}RHS(\Omega(s))ds.
     	\end{equation*}
     	in the closed subspace $Y\subset C\left([\sigma,\sigma+\delta];B_zL_{\xi}^2(m)\right)$ with finite norm
     		\begin{equation*}
     		||\Omega||_{Y}=\sup_{\tau\in[\sigma,\sigma+\delta]}\left(||\Omega(\tau)||_{B_zL_{\xi}^2(m)}+\left(\tau-\sigma\right)^{\frac{1}{2}}\left|\left|\overline{\nabla}\Omega(\tau)\right|\right|_{B_zL_{\xi}^2(m)}\right).
     	\end{equation*}
     	Applying the estimate (\ref{core duhamel estimate}) we see that
     	\begin{equation*}
     		\left|\left|e^{(\tau-\sigma)\mathcal{L}+\left(e^{\tau}-e^{\sigma}\right)\partial_{zz}}\Omega\right|\right|_{Y}\lesssim||\Omega||_{B_zL_{\xi}^2(m)}.
     	\end{equation*}
     	So again matters reduce to proving that the map
     		\begin{equation*}
     			\mathcal{T}:\Omega\mapsto\int_{\sigma}^{\tau}e^{(\tau-s)\mathcal{L}+\left(e^{\tau}-e^s\right)\partial_{zz}}RHS(\Omega(s))ds
     		\end{equation*}
        is s contraction on $Y$ for $\delta>0$ chosen sufficiently small.
        
        In fact,
        \begin{equation*}
        	RHS(\Omega)=-V^{G_{\alpha}}\cdot\nabla_{\xi}\Omega-\Omega\cdot div_{\xi}V^{G_{\alpha}}-V^{\xi}\cdot\nabla_{\xi}G_{\alpha}-G_{\alpha}\cdot\overline{div}V.
        \end{equation*}
        For these four terms, we have the following control
        \begin{align*}
        	& \left|\left|V^{G_{\alpha}}\cdot\nabla_{\xi}\Omega\right|\right|_{B_zL_{\xi}^2(m)}\lesssim(\tau-\sigma)^{-\frac{1}{2}}||\Omega||_{Y},\\
        	&\left|\left|\Omega\cdot div_{\xi}V^{G_{\alpha}}\right|\right|_{B_zL_{\xi}^2(m)}\lesssim\left|\left|\Omega\right|\right|_{Y},\\
        	&\left|\left|G_{\alpha}\cdot\overline{div}V\right|\right|_{B_zL_{\xi}^2(m)}\lesssim\left|\left|\Omega\right|\right|_{Y},\\
        	&\left|\left|V^{\xi}\cdot\nabla_{\xi}G_{\alpha}\right|\right|_{B_zL_{\xi}^2(m)}\lesssim||V||_{B_zL_{\xi}^4}\lesssim||\Omega||_{B_zL_{\xi}^2(m)}\lesssim||\Omega||_{Y}.
        \end{align*}
     	We just proved that
     	\begin{equation*}
     		||RHS(\Omega(s))||_{B_zL_{\xi}^2(m)}\lesssim\left(s-\sigma\right)^{-\frac{1}{2}}||\Omega||_{X}.
     	\end{equation*}
     	Therefore
     	 \begin{align*}
     		||\mathcal{T}\Omega(\tau)||_{B_zL_{\xi}^2(m)}&\lesssim\left(\tau-\sigma\right)^{\frac{1}{2}}||\Omega||_{Y},\\
     		\left(\tau-\sigma\right)^{\frac{1}{2}}||\overline{\nabla}\mathcal{T}\Omega(\tau)||_{B_zL_{\xi}^2(m)}&\lesssim\left(\tau-\sigma\right)^{\frac{1}{2}}||\Omega||_Y.
     	\end{align*}
     	Overall, we find that
     	\begin{equation*}
     		||\mathcal{T}\Omega||_{Y}\lesssim\delta^{\frac{1}{2}}||\Omega||_Y.
     	\end{equation*}
     	so choosing $\delta$ small enough $\mathcal{T}$ is a contraction on $Y$, from which the desired estimate follows.
     \end{proof}
     
     \subsection{Proof of Theorem 3.1}
     From Proposition 3.2 we know that the solution operator $S(\tau,\sigma)$ is well defined and satisfies the estimates (\ref{core estimate 1}), (\ref{decay core estimate}). Thus it remains to prove the estimate (\ref{core estimate 2}).
     
     We first take $\delta>0$ to be the minimum of the $\delta's$ from Propositions 3.7 and 3.8. 
     
     \underline{\textit{Case 1:} $\sigma\le\tau\le\sigma+\delta$}. In this regime, $a(\tau-\sigma)\approx\tau-\sigma$; here the estimate (\ref{core estimate 2}) follows directly from Proposition 3.7. 
     
     \underline{\textit{Case 2:} $\tau>\sigma+\delta$}.
     We first note that in this case $a(\tau-\sigma)\approx 1$. Next we apply Proposition 3.7 on the time interval $[\sigma,\sigma+\frac{\delta}{2}]$ and, writing $	S\left(\sigma+\frac{\delta}{2},\sigma\right)\overline{div}F=\overline{div}R\left(\sigma+\frac{\delta}{2},\sigma\right)F$, we may decompose
     \begin{equation*}
     	S\left(\sigma+\frac{\delta}{2},\sigma\right)\overline{div}=h_1+e^{\frac{\sigma}{2}}\partial_{z}h_2,
     \end{equation*}
     where $h_1$ and $h_2$ are scalars satisfying the estimates
     \begin{equation*}
     	||h_1||_{B_zL_{\xi}^2(m)}+||h_2||_{B_zL_{\xi}^2(m)}\lesssim||F||_{B_zL_{\xi}^p(m)},
     \end{equation*}
     and $\int_{\mathbb{R}^2}h_1d\xi=0$.
     
     For $h_1$ we apply the long time estimate (\ref{decay core estimate}) to obtain
     \begin{equation*}
     	\left|\left|S\left(\tau,\sigma+\frac{\delta}{2}\right)h_1\right|\right|_{B_zL_{\xi}^2(m)}\lesssim e^{-\nu(\tau-\sigma)}||h_1||_{B_zL_{\xi}^2(m)}\lesssim e^{-\nu(\tau-\sigma)}||F||_{B_zL_{\xi}^p(m)}.
     \end{equation*}
     For $h_2$ we instead apply the long time estimate (\ref{core estimate 1}) with the short time estimate in Proposition 3.8 and the fact that $\partial_{z}$ commutes with $S(\tau,\sigma)$ to obtain
     \begin{align*}
     	\left|\left|S\left(\tau,\sigma+\frac{\delta}{2}\right)e^{\frac{\sigma}{2}}\partial_{z}h_2\right|\right|_{B_zL_{\xi}^2(m)}&=e^{\frac{\sigma}{2}}\left|\left|\partial_{z}S\left(\tau,\tau-\frac{\delta}{2}\right)S\left(\tau-\frac{\delta}{2},\sigma+\frac{\delta}{2}\right)h_2\right|\right|_{B_zL_{\xi}^2(m)}\\
     	&\lesssim e^{-\frac{\tau-\sigma}{2}}\left|\left|S\left(\tau-\frac{\delta}{2},\sigma+\frac{\delta}{2}\right)h_2\right|\right|_{B_zL_{\xi}^2(m)}\\
     	&\lesssim e^{-\nu(\tau-\sigma)}||F||_{B_zL_{\xi}^p(m)}.
     \end{align*}
     We complete the proof of Theorem 3.1.
    \section{Linear Estimates for the Background Part}
    In this section, we consider the linearized equation
    \begin{equation}
    	\partial_{t}u+\nabla\cdot(uv^g)=\Delta u, \label{background linear}
    \end{equation}
    where $v^g=\begin{bmatrix}
    	\frac{1}{\sqrt{t}}V^{G_{\alpha}}\left(\frac{x}{\sqrt{t}}\right)\\
    	0
    \end{bmatrix}$. We denote the corresponding evolution operator by $\mathbb{S}(t,s)$.
    
    $\mathbb{S}(t,s)$ has the following estimate:
    
    \bigskip
    
    \textbf{Proposition 4.1.} If $1\le q\le p\le \infty$, then for all $s\in[0,t)$,
    \begin{equation}
    	||\mathbb{S}(t,s)u||_{B_zL_x^p}\lesssim(t-s)^{\frac{1}{p}-1}||u||_{B_zL_x^1}. \label{Prop 4.1}
    \end{equation}
    \begin{proof}
    	Taking $b(t)=e^{-(t-s)|\zeta|^2}\hat{u}(t)$, we obtain the two-dimensional equation:
    	\begin{equation}
    		\partial_{t}b+\nabla_x\cdot\left(b(t)\frac{1}{\sqrt{t}}V^{G_{\alpha}}\left(\frac{x}{\sqrt{t}}\right)\right)=\Delta_xb \label{2D reduction background eq}
    	\end{equation}
    	
    	Now fix the frequency $\zeta$. One may apply Lemma 16 and uniqueness proof of Proposition 6 (i) in \cite{2D PKS} to obtain the following estimate:
    		\begin{equation*}
    		||b(t)||_{L_x^p}\lesssim(t-s)^{\frac{1}{p}-\frac{1}{q}}||b_s||_{L_x^q}.
    	\end{equation*}
    	Switching back to the physical variable, we obtain the desired estimate (\ref{Prop 4.1}).
    \end{proof}

    Now we turn to estimate $\mathbb{S}(t,s)div$. Note that 
    \begin{equation}
    	\mathbb{S}(t,s)div=divK(t,s) \label{background commute}
    \end{equation}
    where $K(t,s)$ is the evolution operator of the following equation:
    \begin{equation}
    	\partial_{t}F+v^gdivF=\Delta F \label{background K(t,s) eq}
    \end{equation}
    We need to estimate $K(t,s)$:
    
    \bigskip
    
    \textbf{Proposition 4.2.} For all $s\in(0,t)$, we have
    \begin{equation}
    	||K(t,s)f||_{B_zL_x^1}\lesssim_{\alpha}||f||_{B_zL_x^1}.
    \end{equation}
    
    \begin{proof}
    	Taking the $z$-direction Fourier transform in (\ref{background K(t,s) eq}), we obtain the equation
    	\begin{equation*}
    		\partial_{t}\hat{F}+v^g\left(div_x\hat{F^x}+i\zeta\hat{F^z}\right)=\left(\Delta_x-|\zeta|^2\right)\hat{F}.
    	\end{equation*}
    	
    	We may then decompose this into a system of two equations:
    		\begin{equation*}
    		\begin{cases}
    			\partial_{t}\hat{F}^x+\left(v^g\right)^x\left(div_x\hat{F}^x+i\zeta\hat{F}^z\right)=\left(\Delta_x-|\zeta|^2\right)\hat{F}^x,\\
    			\partial_{t}\hat{F}^z=\left(\Delta_x-|\zeta|^2\right)\hat{F}^z.
    		\end{cases}
    	\end{equation*}
    	
    	Next we switch to self-similar variables, letting
    	\begin{equation*}
    		H(\tau,\xi,\zeta)=e^{\frac{\tau}{2}}\hat{F}\left(e^{\tau},e^{\frac{\tau}{2}}\xi,\zeta\right).
    	\end{equation*}
    	to obtain the system
    	\begin{equation*}
    		\begin{cases}
    			\partial_{\tau}H^x+V^{G_{\alpha}}(\xi)div_{\xi}H^x+i\zeta e^{\frac{\tau}{2}}V^{G_{\alpha}}(\xi)H^z=\left(\mathcal{L}-\frac{1}{2}-e^{\tau}|\zeta|^2\right)H^x,\\
    			\partial_{\tau}H^z=\left(\mathcal{L}-\frac{1}{2}-e^{\tau}|\zeta|^2\right)H^z.
    		\end{cases}
    	\end{equation*}
    	\begin{itemize}
    		\item $H^z$ part: Using the semigroup estimate of $e^{\tau\mathcal{L}}$ on $L^1$,
    		\begin{equation*}
    			||H^z(\tau)||_{L_{\xi}^1}\lesssim e^{-\frac{1}{2}(\tau-\sigma)}||H^z(\sigma)||_{L_{\xi}^1}.
    		\end{equation*}
    		\item $H^x$ part: Denote
    		\begin{equation*}
    			\Xi_{\alpha}:=V^{G_{\alpha}}(\xi)div_{\xi}.
    		\end{equation*}
    		A contraction mapping argument similar to Proposition 3.7 shows that there exists $0<\delta=\delta(\alpha)<<1$ such that for all $\tau\in[0,\delta]$,
    		\begin{equation*}
    			||e^{\tau(\mathcal{L}-\Xi_{\alpha})}||_{L_{\xi}^1}\lesssim_{\alpha}||f||_{L_{\xi}^1}.
    		\end{equation*}
    		Obviously this estimate also holds for all $\tau>0$. Now write the integral equation for $H^x$:
    		\begin{equation*}
    			H^x(\tau)=e^{-(e^{\tau}-e^{\sigma})|\zeta|^2+(\tau-\sigma)(\mathcal{L}-\Xi_{\alpha}-\frac{1}{2})}H^x(\sigma)-\int_{\sigma}^{\tau}e^{-(e^{\tau}-e^{s})|\zeta|^2+(\tau-s)(\mathcal{L}-\Xi_{\alpha}-\frac{1}{2})}i\zeta e^{\frac{s}{2}}V^{G_{\alpha}}(\xi)H^z(s)ds.
    		\end{equation*}
    		Taking the $L_{\xi}^1$ norm on both sides and using the estimate for $H^z$ part, we obtain
    		\begin{align*}
    			||H^x(\tau)||_{L_{\xi}^1}&\lesssim e^{-\frac{1}{2}(\tau-\sigma)}||H^x(\sigma)||_{L_{\xi}^1}+\int_{\sigma}^{\tau}e^{-(e^{\tau}-e^s)|\zeta|^2-\frac{1}{2}(\tau-s)}|\zeta|e^{\frac{s}{2}}||H^z(s)||_{L_{\xi}^1}ds\\
    			&\lesssim e^{-\frac{1}{2}(\tau-\sigma)}||H^x(\sigma)||_{L_{\xi}^1}+\int_{\sigma}^{\tau}e^{-(e^{\tau}-e^s)|\zeta|^2-\frac{1}{2}(\tau-s)}|\zeta|ds\cdot e^{\frac{\sigma}{2}}||H^z(\sigma)||_{L_{\xi}^1}\\
    			&= e^{-\frac{1}{2}(\tau-\sigma)}||H^x(\sigma)||_{L_{\xi}^1}+\int_{\sigma}^{\tau}e^{-(e^{\tau}-e^s)|\zeta|^2-\frac{1}{2}(\tau-s)}|\zeta|\sqrt{e^{\tau}-e^s}\dfrac{1}{\sqrt{e^{\tau}-e^s}}ds\cdot e^{\frac{\sigma}{2}}||H^z(\sigma)||_{L_{\xi}^1}\\
    			&\lesssim e^{-\frac{1}{2}(\tau-\sigma)}||H^x(\sigma)||_{L_{\xi}^1}+\int_{\sigma}^{\tau}\dfrac{1}{\sqrt{e^{\tau-s}-1}}ds\cdot e^{-\frac{1}{2}(\tau-\sigma)}||H^z(\sigma)||_{L_{\xi}^1}\\
    			&\lesssim e^{-\frac{1}{2}(\tau-\sigma)}(||H^x(\sigma)||_{L_{\xi}^1}+||H^z(\sigma)||_{L_{\xi}^1})
    		\end{align*}
    	\end{itemize}
    	Combining the above estimates and returning to the physical variables, we obtain the desired estimate for $K(t,s)$.
    \end{proof}
    We also need short time estimates for $\mathbb{S}(t,s)div$.
    
    \bigskip
    
    \textbf{Proposition 4.3.} There exists $0<\delta=\delta(\alpha)<<1$ so that for all $0<s<t\le(1+\delta)s$ and all 3-vector $F\in B_zL_x^1$, we have the estimate
    \begin{equation}
    	||\mathbb{S}(t,s)divF||_{B_zL_x^1}\lesssim\left(t-s\right)^{-\frac{1}{2}}||F||_{B_zL_x^1}.
    \end{equation}
    
    \begin{proof}
    	The proof is similar to the short time estimates in Section 3. Write (\ref{background K(t,s) eq}) in integral form
    	\begin{equation*}
    		F(t)=e^{(t-s)\Delta}F(s)-\int_{s}^{t}e^{(t-\sigma)\Delta}v^gdivF(\sigma)d\sigma.
    	\end{equation*}
    	Define
    		\begin{align*}
    		&\mathcal{T}F=e^{(t-s)\Delta}F(s)-\int_{s}^{t}e^{(t-\sigma)\Delta}v^gdivF(\sigma)d\sigma,\\
    		&Z=\{F\in B_zL_x^1:||F||_Z<+\infty\},\\
    		&||F||_Z=\sup_{t\in[s,(1+\delta)s]}\left(||F||_{B_zL_x^1}+(t-s)^{\frac{1}{2}}||divF||_{B_zL_x^1}\right).
    	\end{align*}
    	Using the heat kernel estimates
    		\begin{align*}
    		&||e^{t\Delta}f||_{B_zL_x^1}\lesssim||f||_{B_zL_x^1},\\
    		&||e^{t\Delta}divf||_{B_zL_x^1}\lesssim t^{-\frac{1}{2}}||f||_{B_zL_x^1},
    	\end{align*}
    	We obtain
    		\begin{align*}
    		||TF(t)||_{B_zL_x^1}&\lesssim||F(s)||_{B_zL_x^1}+\int_{s}^{t}\left(t-\sigma\right)^{-\frac{1}{2}}||F\otimes v^g||_{B_zL_x^1}d\sigma+\int_{s}^{t}||\left(F\cdot\nabla\right)v^g||_{B_zL_x^1}d\sigma\\
    		&\lesssim||F(s)||_{B_zL_x^1}+\int_{s}^{t}\left(t-\sigma\right)^{-\frac{1}{2}}\sigma^{-\frac{1}{2}}||F||_Zd\sigma+\log(1+\delta)||F||_Z\\
    		&\lesssim||F||_Z\left(1+\sqrt{\delta}+\log(1+\delta)\right),
    	\end{align*}
    	\begin{align*}
    		\left(t-s\right)^{\frac{1}{2}}||divTF||_{B_zL_x^1}&\lesssim||F(s)||_{B_zL_x^1}+\left(t-s\right)^{\frac{1}{2}}\int_{s}^{t}\left(t-\sigma\right)^{-\frac{1}{2}}||v^gdivF(\sigma)||_{B_zL_x^1}d\sigma\\
    		&\lesssim||F(s)||_{B_zL_x^1}+\left(t-s\right)^{\frac{1}{2}}\int_{s}^{t}\left(t-\sigma\right)^{-\frac{1}{2}}\sigma^{-\frac{1}{2}}\left(\sigma-s\right)^{-\frac{1}{2}}d\sigma||F||_Z\\
    		&\lesssim||F||_{Z}\left(1+\sqrt{\delta}\right).
    	\end{align*}
    	Hence $\mathcal{T}$ is a contraction on $Z$. Banach's fixed point theorem tells us that there exists a unique solution of (\ref{background K(t,s) eq}) in $Z$ satisfying
    		\begin{equation*}
    		\left(t-s\right)^{\frac{1}{2}}||divF||_{B_zL_x^1}\lesssim||F(s)||_{B_zL_x^1}+\sqrt{\delta}||F||_Z.
    	\end{equation*}
    	Taking the supremum of $t$, we get
    	\begin{equation*}
    		\left(1-\sqrt{\delta}\right)||F||_Z\lesssim||F(s)||_{B_zL_x^1}.
    	\end{equation*}
    	Therefore
    	\begin{equation*}
    		||\mathbb{S}(t,s)divF||_{B_zL_x^1}\lesssim\left(t-s\right)^{-\frac{1}{2}}||F||_{B_zL_x^1}.
    	\end{equation*}
    \end{proof}
    Finally we establish the long time estimate for $\mathbb{S}(t,s)div$:
    
    \bigskip
    
    \textbf{Theorem 4.4.} Let $0<s<t\le T$, $F$ is a 3-vector in $B_zL_x^1$, then for all $p\in[1,\infty]$ we have the estimate
    \begin{equation}
    	||\mathbb{S}(t,s)divF||_{B_zL_x^p}\lesssim\left(t-s\right)^{-\left(\frac{3}{2}-\frac{1}{p}\right)}||F||_{B_zL_x^1}, \label{background main estimate div}
    \end{equation}
    The implicit constant depends on $p,\alpha$.
    \begin{proof}
    	First notice that if $p=1$ case holds, then for general $p$ we have
    		\begin{align*}
    		||\mathbb{S}(t,s)divF||_{B_zL_x^p}&=\left|\left|\mathbb{S}\left(t,\frac{t+s}{2}\right)\mathbb{S}\left(\frac{t+s}{2},s\right)divF\right|\right|_{B_zL_x^p}\\
    		&\lesssim(t-s)^{-\left(1-\frac{1}{p}\right)}\left|\left|\mathbb{S}\left(\frac{t+s}{2},s\right)divF\right|\right|_{B_zL_x^1}\\
    		&\lesssim\left(t-s\right)^{-\left(\frac{3}{2}-\frac{1}{p}\right)}||F||_{B_zL_x^1}.
    	\end{align*}
    	Hence it suffices to prove the case for $p=1$. Take the $\delta$ in Proposition 4.3.
    	
    	\underline{\textit{Case 1:} $s<t\le(1+\delta)s$}. The estimate can be derived directly from Proposition 4.3.
    	
    	\bigskip
    	
    	\underline{\textit{Case 2:} $(1+\delta)s<t\le T$}. In this case we apply (\ref{background commute}) near the end point $T$:
    	\begin{align*}
    		||\mathbb{S}(t,s)divF||_{B_zL_x^1}&\lesssim\left|\left|\mathbb{S}\left(t,\frac{2}{2+\delta}t\right)divK\left(\frac{2}{2+\delta}t,s\right)F\right|\right|_{B_zL_x^1}\\
    		&\lesssim t^{-\frac{1}{2}}\left|\left|K\left(\frac{2}{2+\delta}t,s\right)F\right|\right|_{B_zL_x^1}\\
    		&\lesssim
    		\left(t-s\right)^{-\frac{1}{2}}||F||_{B_zL_x^1}.
    	\end{align*}
    \end{proof}
    \section{Nonlinear Estimates}
    This section is devoted to the proof of Theorem 2.2. We start by splitting the cell density
    \begin{equation*}
    	u(t,x,z)=\frac{1}{t}G_{\alpha}\left(\frac{x}{\sqrt{t}}\right)+u^c(t,x,z)+u^b(t,x,z)
    \end{equation*}
    with the corresponding velocity field given by the Biot-Savart law
    \begin{equation*}
    	v(t,x,z)=v^g(t,x,z)+v^c(t,x,z)+v^b(t,x,z)
    \end{equation*}
    
    As usual, we capitalize $u$ and $v$ in self-similar variables: for $*=b\ \text{or}\ c$,
    \begin{equation*}
    	U^*(\tau,\xi,z)=tu^*(t,x,z),\ \ \ V^*(\tau,\xi,z)=\sqrt{t}v^*(t,x,z).
    \end{equation*}
    
    Recall that $S$ is the semigroup associated to $\partial_{\tau}U+\nabla_{\xi}\cdot(U V^{G_{\alpha}})+\overline{\nabla}\cdot(G_{\alpha}V)=(\mathcal{L}+e^{\tau}\partial_{zz})U$, while $\mathbb{S}$ is the semigroup associated to $	\partial_{t}u+\nabla\cdot(uv^g)=\Delta u$.
    
    Duhamel's formula then formally gives
    \begin{equation*}
    	u^b(t)=\mathbb{S}(t,0)\mu^b-\int_{0}^{t}\mathbb{S}(t,s)\nabla\cdot\left(u^b{v}^c\right)ds-\int_{0}^{t}\mathbb{S}(t,s)\nabla\cdot\left(u^bv^b\right)ds,
    \end{equation*}
    \begin{equation*}
    	U^c(\tau)=-\int_{-\infty}^{\tau}S(\tau,\sigma)\overline{\nabla}\cdot\left(G_{\alpha}V^b\right)d\sigma-\int_{-\infty}^{\tau}S(\tau,\sigma)\overline{\nabla}\cdot\left(U^c{V}^b\right)d\sigma-\int_{-\infty}^{\tau}S(\tau,\sigma)\overline{\nabla}\cdot\left(U^cV^c\right)d\sigma.
    \end{equation*}
    Writing $\mathcal{Q}$ for the above right-hand side, we are looking for a solution of the equation
    \begin{equation*}
    	(u^b,U^c)=\left(\mathcal{Q}^b\left(u^b,U^c\right),\mathcal{Q}^c\left(u^b,U^c\right)\right)=\mathcal{Q}\left(u^b,U^c\right).
    \end{equation*}
    Taking $m>2$, we will solve this fixed point problem by applying the contraction mapping principle in the following ball, for constants $M,D\ge 1$ determined by the proof below,
    \begin{equation*}
    	B_{\epsilon,T}=\left\{(u^b,U^c)\bigg|\left|\left|(u^b,U^c)\right|\right|_{X}:=M\sup_{t\in(0,T)}t^{\frac{1}{4}}\left|\left|u^b(t)\right|\right|_{B_zL_x^{\frac{4}{3}}}+\sup_{\tau\in(-\infty,\log T)}\left|\left|U^c(\tau)\right|\right|_{B_zL_{\xi}^2(m)}\le D\epsilon\right\},
    \end{equation*}
    
    To prove Theorem 2.2 we will verify that, whenever the data satisfies $||\mu^b||_{B_zL_x^1}\le\epsilon$, the map $\mathcal{Q}:B_{\epsilon,T}\rightarrow B_{\epsilon,T}$ is a contraction for any $0<\epsilon\le\epsilon_0$, any $T>0$ and an appropriate choice of the constants $M,D$ and $\epsilon_0$.
    
    \bigskip
    
    \underline{Bound for the core}. We abbreviate the three summands in the definition of $\mathcal{Q}^c$ by
    	\begin{equation*}
    	Q^c(u^b,U^c)=L^c+N_0^c+N_1^c.
    \end{equation*}
    By (\ref{core estimate 2}) and the rapid decay of $G_{\alpha}$,
    \begin{equation*}
    	\left|\left|L^c(\tau)\right|\right|_{B_zL_{\xi}^2(m)}\lesssim\int_{-\infty}^{\tau}\dfrac{e^{-\nu(\tau-\sigma)}}{\left(a(\tau-\sigma)\right)^{\frac{1}{2}}}\left|\left|G_{\alpha}V^b(\sigma)\right|\right|_{B_zL_{\xi}^2(m)}d\sigma\lesssim\int_{-\infty}^{\tau}\dfrac{e^{-\nu(\tau-\sigma)}}{\left(a(\tau-\sigma)\right)^{\frac{1}{2}}}\left|\left|V^b(\sigma)\right|\right|_{B_zL_{\xi}^4}d\sigma
    \end{equation*}
    By scaling and the Biot-Savart bounds in Lemma 3.5 (iii)
    \begin{equation*}
    	\left|\left|V^b(\tau)\right|\right|_{B_zL_{\xi}^4}=e^{\frac{\tau}{4}}\left|\left|v^b(e^{\tau})\right|\right|_{B_zL_x^4}\lesssim t^{\frac{1}{4}}\left|\left|u^b(t)\right|\right|_{B_zL_x^{\frac{4}{3}}}.
    \end{equation*}
    Hence, with an implicit constant independent of $T$ and $\epsilon$,
    \begin{equation*}
    ||L^c(\tau)||_{B_zL_{\xi}^2(m)}\lesssim\frac{1}{M}\left|\left|(u^b,U^c)\right|\right|_X\int_{-\infty}^{\tau}\dfrac{e^{-\nu(\tau-\sigma)}}{\left(a(\tau-\sigma)\right)^{\frac{1}{2}}}d\sigma\lesssim\frac{1}{M}\left|\left|(u^b,U^c)\right|\right|_X.
    \end{equation*}
    Using again (\ref{core estimate 2}), we have
    	\begin{align*}
    	\left|\left|N_0^c\right|\right|_{B_zL_{\xi}^2(m)}&\lesssim\int_{-\infty}^{\tau}\dfrac{e^{-\nu(\tau-\sigma)}}{(a(\tau-\sigma))^{\frac{3}{4}}}\left|\left|U^c{V}^b\right|\right|_{B_zL_{\xi}^{4/3}(m)}d\sigma\\	
    	&\lesssim\int_{-\infty}^{\tau}\dfrac{e^{-\nu(\tau-\sigma)}}{(a(\tau-\sigma))^{\frac{3}{4}}}\left|\left|U^c\right|\right|_{B_zL_{\xi}^2(m)}\left|\left|{V}^b\right|\right|_{B_zL_{\xi}^4}d\sigma\\
    	&\lesssim\frac{1}{M}\left|\left|(u^b,U^c)\right|\right|_X^2\int_{-\infty}^{\tau}\dfrac{e^{-\nu(\tau-\sigma)}}{(a(\tau-\sigma))^{\frac{3}{4}}}d\sigma\\
    	&\lesssim\left|\left|(u^b,U^c)\right|\right|_X^2.
    \end{align*}
    The same proof applies to the nonlinear term $N_1^c$:
    	\begin{equation*}
    	\left|\left|N_1^c\right|\right|_{B_zL_{\xi}^2(m)}\lesssim\left|\left|(u^b,U^c)\right|\right|_X^2.
    \end{equation*}
    
    \bigskip
    
    \underline{Bound for the background}. We abbreviate the three summands in the definition of $\mathcal{Q}^b$ by
    \begin{equation*}
    	Q^b\left(u^b,U^c\right)=L^b+N_0^b+N_1^b.
    \end{equation*}
   By (\ref{background main estimate div}),
    \begin{align*}
    	t^{\frac{1}{4}}\left|\left|N_0^b(t)\right|\right|_{B_zL_x^{\frac{4}{3}}}&\lesssim t^{\frac{1}{4}}\int_{0}^{t}\dfrac{1}{\left(t-s\right)^{\frac{3}{4}}}\left|\left|u^b{v}^c\right|\right|_{B_zL_x^1}ds\\
    	&\lesssim t^{\frac{1}{4}}\int_{0}^{t}\dfrac{1}{\left(t-s\right)^{\frac{3}{4}}}\left|\left|u^b\right|\right|_{B_zL_x^{\frac{4}{3}}}\left|\left|{v}^c\right|\right|_{B_zL_x^4}ds\\
    	&\lesssim t^{\frac{1}{4}}\frac{1}{M}\left|\left|\left(u^b,U^c\right)\right|\right|_X^2\int_{0}^{t}\dfrac{1}{\left(t-s\right)^{\frac{3}{4}}s^{\frac{1}{2}}}ds\\
    	&\lesssim\frac{1}{M}\left|\left|\left(u^b,U^c\right)\right|\right|_X^2. 
    \end{align*}
   where we used scaling and Lemma 3.5 (iii) 
    \begin{equation*}
    	||v^c(t)||_{B_zL_x^4}\lesssim||u^c(t)||_{B_zL_x^{\frac{4}{3}}}=t^{-\frac{1}{4}}||U^c(\tau)||_{B_zL_{\xi}^{\frac{4}{3}}}\lesssim t^{-\frac{1}{4}}||U^c(\tau)||_{B_zL_{\xi}^2(m)}\lesssim t^{-\frac{1}{4}}||(u^b,U^c)||_X.
    \end{equation*}
    Similarly,
    	\begin{equation*}
    	t^{\frac{1}{4}}\left|\left|N_1^b(t)\right|\right|_{B_zL_x^{\frac{4}{3}}}\lesssim\frac{1}{M^2}\left|\left|(u^b,U^c)\right|\right|_X^2.
    \end{equation*}
    
    \bigskip
    
    Combining Proposition 4.1 and the above estimates, we obtain
    \begin{align*}
    	\left|\left|Q(u^b,U^c)\right|\right|_X&\lesssim M\left|\left|\mu^b\right|\right|_{B_zL_x^1}+\left|\left|(u^b,U^c)\right|\right|_X^2\left(1+\frac{1}{M}\right)+\left|\left|(u^b,U^c)\right|\right|_X\frac{1}{M}\\
    	&\lesssim M\left|\left|\mu^b\right|\right|_{B_zL_x^1}+\frac{D\epsilon}{M}+\frac{M+1}{M}D^2\epsilon^2.
    \end{align*}
    Let
    \begin{equation*}
    	\left|\left|\mu^b\right|\right|_{B_zL_x^1}\le\epsilon,\ \ \ M+\frac{D}{M}<D,\ \ \ \epsilon\le\dfrac{D-M-\frac{D}{M}}{\frac{M+1}{M}D^2},\ \ \ \epsilon\le\frac{\beta-\frac{1}{M}}{(1+\frac{1}{M})D}\ \ \left(\beta\in(0,1)\right).
    \end{equation*}
    Then $\mathcal{Q}:B_{\epsilon,T}\rightarrow B_{\epsilon,T}$ is a contraction, which proves the main theorem.
    
    \bigskip
    
    \underline{Mild solution of PKS}. Finally, we verify that the solution constructed by the contraction mapping argument above does indeed satisfy Definition 2.1 (the proof is identical to that in \cite{filament}).  For this we need two things: $(a)\ u\in C_w([0,T];M^{\frac{3}{2}})$, and in particular $\lim_{t\searrow 0}u(t)=\alpha\delta_{x=0}+\mu^b$ in the sense of $\mathscr{S}'$; $(b)$ that the mild form of the equations (\ref{duhamel formula}) is satisfied in the sense of tempered distributions.
    
    Using Proposition 4.1, the contraction mapping argument above and the embedding $B_zL_x^1\subset M^{\frac{3}{2}}$, one can obtain the a priori estimate 
    \begin{equation*}
    	\sup_{0<t<T}||u(t)||_{M^{\frac{3}{2}}}\lesssim 1
    \end{equation*}
    and classical parabolic regularity ensures that $u\in C((0,T],M^{\frac{3}{2}})$.
    
    It remains to verify that $\lim_{t\searrow 0}u(t)=\alpha\delta_{x=0}+\mu^b$ in the sense of tempered distributions. The term involving $G_{\alpha}$ converges to $\alpha\delta_{x=0}$ weak* (hence in the sense of tempered distributions). By an inspection of the contraction mapping argument above, we see that the solution constructed satisfies
    \begin{equation*}
    	\lim\limits_{t\searrow 0}||u^c(t)||_{B_zL_x^1}=0.
    \end{equation*}
     Appealing to the Duhamel formula for $u^b$, for any $\phi\in\mathscr{S}$, we may apply (\ref{background commute}) to yield 
    \begin{equation*}
    	\int\phi\cdot\mathbb{S}(t,s)\nabla\cdot(u^bv^b)dx
    	=-\int\nabla\phi\cdot K(t,s)(u^bv^b)dx
    	\lesssim||u^b(s)||_{B_zL_x^4}||\omega^b(s)||_{B_zL_x^{\frac{4}{3}}}\lesssim\frac{1}{s^{\frac{1}{2}}}.
    \end{equation*}
    and similarly for the other nonlinear term in the $u^b$ equation. Hence, by the weak* continuity of $\mathbb{S}(t,s)$, $u^b(s)$ converges to $\mu^b$ in $\mathscr{S}'$ and we have proved that $u\in C_w([0,T];M^{\frac{3}{2}})$.
    
    Next we verify that $u(t)$ satisfies (\ref{duhamel formula}). By parabolic regularity, $u$ is a classical solution of the 3D Patlak-Keller-Segel system for $t>0$. Hence, for all $0<s<t$,
    \begin{equation*}
    	u(t)=e^{(t-s)\Delta}u(s)-\int_{s}^{t}e^{(t-t')\Delta}\text{div}(u(t')\nabla c(t'))dt'.
    \end{equation*}
   
    Due to the self-adjointness of the heat semigroup in $L^2$, we can pass to the limit in the first term due to the continuity of $u$ in $\mathscr{S}'$: $\lim_{s\searrow 0} e^{(t-s)\Delta}u(s)=e^{t\Delta}u(0)$. Furthermore, we may pass $s\searrow 0$ also in the nonlinear term: indeed, using heat semigroup estimates on Morrey spaces (see \cite{mild solution morrey} Proposition 3.2), we have
    \begin{align*}
    	\int_{0}^{t}\left|\left|e^{(t-t')\Delta}\text{div}(u(t')\nabla c(t'))\right|\right|_{M^{\frac{3}{2}}}dt'&\lesssim\int_{0}^{t}\dfrac{1}{(t-t')^{\frac{1}{2}}}||u(t')\nabla c(t')||_{M^{\frac{3}{2}}}dt'\\
    	&\lesssim\left(1+||(u^b,U^c)||_{X}+||(u^b,U^c)||^2_X\right)\int_{0}^{t}\dfrac{1}{(t-t')^{\frac{1}{2}}(t')^{\frac{1}{2}}}dt'.
    \end{align*}
     \bigskip
    
    \textbf{Acknowledgements.} The author would like to deeply thank Prof. Lifeng Zhao and Dr. Dengjun Guo for fruitful discussions during the pursue of master's degree at University of Science and Technology of China from 2021 to 2023.

\end{document}